\documentclass[11pt]{article}
\usepackage[dvips]{graphics}
\usepackage{multirow}
\usepackage{array}
\usepackage{epsfig,rotating}
\usepackage{amssymb,amsmath}
\usepackage{latexsym}
\usepackage[usenames]{color}

\numberwithin{equation}{section}

\newcommand{\be}{\begin{equation}}
\newcommand{\ee}{\end{equation}}
\newcommand{\beaa}{\begin{eqnarray*}}
\newcommand{\eeaa}{\end{eqnarray*}}
\newcommand{\bea}{\begin{eqnarray}}
\newcommand{\eea}{\end{eqnarray}}

\newcommand{\bei}{\begin{itemize}}
\newcommand{\eei}{\end{itemize}}

\newcommand{\bd}{\mathbf}

\newcommand{\bz}{\mathbf{z}}

\newtheorem{theorem}{\noindent T{\footnotesize HEOREM}}

\newtheorem{lemma}{\noindent L{\footnotesize EMMA}}[section]

\begin{document}
%\nocite{*}

\title{\bf Limiting Spectral Radii of Circular Unitary Matrices under Light Truncation}
\author{Yu Miao$^{1}$
and  Yongcheng Qi$^{2,*}$\\
Henan Normal University and
University of Minnesota Duluth
}

\date{}
\maketitle

\footnotetext[1]{College of Mathematics and Information Science, Henan Normal University, Henan Province 453007,  China. Email: yumiao728@gmail.com
 %\newline  \indent \  \
}

\footnotetext[2]{Department of Mathematics and Statistics, University of Minnesota Duluth, MN 55812, USA. Email: yqi@d.umn.edu (corresponding author)
%\newline  \indent \ \
%The research of Yongcheng Qi was supported in part by NSF Grant DMS-1916014
}

\begin{abstract}
\noindent  Consider a truncated circular unitary matrix which is a
$p_n$ by $p_n$ submatrix of an $n$ by $n$ circular unitary matrix
after deleting the last $n-p_n$ columns and rows.  Jiang and
Qi~\cite{JiangQi2017} and Gui and Qi~\cite{GQ2018} study the
limiting distributions of the maximum absolute value of the
eigenvalues (known as spectral radius) of the truncated matrix. Some
limiting distributions for the spectral radius for the truncated
circular unitary matrix have been obtained under the following
conditions: (1). $p_n/n$ is bounded away from $0$ and $1$; (2).
$p_n\to\infty$ and $p_n/n\to 0$ as $n\to\infty$; (3). $(n-p_n)/n\to
0$ and $(n-p_n)/(\log n)^3\to\infty$ as $n\to\infty$; (4).
$n-p_n\to\infty$ and $(n-p_n)/\log n\to 0$ as $n\to\infty$; and (5).
$n-p_n=k\ge 1$ is a fixed integer. The spectral radius converges in
distribution to the Gumbel distribution under the first four
conditions and to a reversed Weibull distribution under the fifth
condition.  Apparently, the conditions above do not cover the case
when $n-p_n$ is of order between $\log n$ and $(\log n)^3$. In this
paper, we prove that the spectral radius converges in distribution
to the Gumbel distribution as well in this case, as conjectured by
Gui and Qi~\cite{GQ2018}.
\end{abstract}

%{\red [see, e.g., Theorem 2.8.1
% on page 110 from Andrews, Askey and Roy (1999)]}\\\\


%\noindent\textbf{AMS 2000 Subject Classification: \/}  15A52, 60F99, 60G55, 60G70. \\


\newpage

\section{Introduction}\label{intro}

The study of large random matrices can date back to nearly a century
ago, and one example is Wishart's~\cite{Wishart} work on statistical
properties for large covariance matrices. The theory of random
matrices has been rapidly developed in last few decades and has
found applications in heavy-nuclei atoms (Wigner~\cite{Wigner}),
number theory (Mezzadri and Snaith~\cite{MS}), quantum mechanics
(Mehta~\cite{Mehta}), condensed matter physics
(Forrester~\cite{For}), wireless communications (Couillet and
Debbah~\cite{CD}), to just mention a few.

Statistical properties of large random matrices including their
empirical spectral distributions and spectral radii (the largest
eigenvalues) are of particular interest in the study. For the three
Hermitian matrices including  Gaussian orthogonal ensemble, Gaussian
unitary ensemble and Gaussian symplectic ensemble, Tracy and Widom
\cite{Tracy94, Tracy96} show that their spectral radii converge in
distribution to Tracy-Widom laws. For more consequent applications
of Tracy-Widom laws, see, e.g., Baik et al.~\cite{Baik}, Tracy and
Widom~\cite{TW02}, Johansson~\cite{Johansson07},
Johnstone~\cite{John01, John08}, and Jiang~\cite{Jiang09}. For a
non-Hermitian matrix,  the largest absolute value of its eigenvalues
is refereed to as the spectral radius. The spectral radii for the
real, complex and symplectic Ginibre ensembles are explored by
Rider~\cite{Rider, Rider09_23} and Rider and Sinclair~\cite{RS}, and
their limiting distributions are usually the Gumbel distributions
instead of the Tracy-Widow laws.

In this paper, we are interested in the truncation of the circular
unitary ensemble. The circular unitary ensemble is a random square
matrix with Haar measure on the unitary group, and it is also called
Haar-invariant unitary matrix. Truncations of large Haar unitary
matrices are employed to describe quantum systems with absorbing
boundaries (Casati et al.~\cite{Casati1999}) and have applications
in optical and semiconductor superlattices (Gl\"uck et
al.~\cite{Gluck2002}) and quantum conductance
(Forrester~\cite{Forrester2006}), among many others. More references
on applications can be found in Dong et al.~\cite{Dong}.

Let $\bd{U}$ be an $n\times n$ circular unitary matrix. The $n$
eigenvalues of the circular unitary matrix $\bd{U}$ are distributed
over $\{z\in\mathcal{C}: |z|=1\}$ , where $\mathcal{C}$ is the
complex plane, and their joint density function is given by
 \[
\frac{1}{n!(2\pi)^n}\cdot\prod_{1\leq j<k\leq n}|z_j - z_k|^2;
\]
see, e.g., Hiai and Petz~\cite{HF}.    For integer $p$ with $1\le p<n$,  partition $\bd{U}$ as follows
\begin{eqnarray*}
\bd{U}=\begin{pmatrix}
\bd{A}\ \ \bd{C}^*\\
\bd{B}\ \ \bd{D}
\end{pmatrix}
\end{eqnarray*}
where $\bd{A}$, as a truncation of $\bd{U}$, is a $p\times p$
submatrix. Let $\bz_1, \cdots, \bz_p$ be the $p$ eigenvalues of
$\bd{A}$. According to \.Zyczkowski and Sommers~\cite{Zski}, their
density function is
\begin{equation}\label{Good_forever}
C\cdot\prod_{1\leq j<k\leq p}|z_j - z_k|^2\prod_{j=1}^{p}(1-|z_j|^2)^{n-p-1}I(|z_j|<1)
\end{equation}
where $C$ is a constant, depending on both $n$ and $p$ such that the above function is a probability density.

In this paper we assume that $p=p_n$ depends on $n$ and $\lim_{n\to\infty}p_n=\infty$.

Set $c=\lim_{n\to\infty}(p_n/n)$.  \.{Z}yczkowski and
Sommers~\cite{Zski} prove that the empirical spectral distribution
of $\bz_i$'s converges to the distribution with density proportional
to $\frac{1}{(1-|z|^2)^2}$ for $|z|\leq c$ if $c\in (0,1)$. Dong et
al.~\cite{Dong} show that the empirical spectral distribution goes
to the circular law and the arc law as $c=0$ and $c=1$,
respectively. For more work, see also Diaconis and
Evans~\cite{Evans} and Jiang~\cite{Jiang09, Jiang10}.

Two recent papers by Jiang and Qi~\cite{JiangQi2017} and Gui and
Qi~\cite{GQ2018} study the limiting distributions of the spectral
radius $\max_{1\le j\le p}|\bz_j|$ for the truncated circular
unitary ensemble.  Jiang and Qi~\cite{JiangQi2017} have proved that
the spectral radius $\max_{1\le j\le p}|\bz_j|$ converges to the
Gumbel distribution when the ratio $p_n/n$ is bounded away from $0$
and $1$.  Gui and Qi~\cite{GQ2018} further consider the case when
the limit of $p_n/n$ is $0$ or $1$. Since $\bd{A}$ is obtained by
deleting last $n-p_n$ rows and columns from $\bd{U}$, we call the
truncation is light if $\lim_{n\to\infty}(n-p_n)/n=0$, otherwise the
truncation is heavy if $\liminf_{n\to\infty}(n-p_n)/n>0$.  The main
results obtained by Jiang and Qi~\cite{JiangQi2017} and Gui and
Qi~\cite{GQ2018} are summarized in section~\ref{main}.

Obvious, from Jiang and Qi~\cite{JiangQi2017} and Gui and
Qi~\cite{GQ2018} we observe that the limiting distribution for the
spectral radius $\max_{1\le j\le p}|\bz_j|$ depends on the
truncation parameter $n-p_n$.  We are interested in investigating
how the limiting distribution of the spectral radius changes when
the the truncation parameter runs over the range $1\le n-p_n<n$
under constrain that $\lim_{n\to\infty}p_n=\infty$. The limiting
distribution for the spectral radius $\max_{1\le j\le p}|\bz_j|$
remains unknown when the truncation parameter $n-p_n$ is of order
between $\log n$ and $(\log n)^3$. Gui and Qi~\cite{GQ2018}
conjecture that $\max_{1\leq j \leq p_n}|\bz_j|$, after properly
normalized, converges in distribution to the Gumbel distribution in
this case. In this paper, we will show that this conjecture is true.
This paper together with Jiang and Qi~\cite{JiangQi2017} and Gui and
Qi~\cite{GQ2018} will put an end to the study of the limiting
spectral radius for the truncated circular unitary ensemble. It is
worth noting that the key approaches for the proofs in Jiang and
Qi~\cite{JiangQi2017} and Gui and Qi~\cite{GQ2018} are no longer
applicable in the aforementioned regime, and therefore, we have to
use a totally different approach in this paper. More details will be
provided in {\bf Remark} 2 in section~\ref{main}.

The rest of the paper is organized as follows.  The main result in
this paper is given in section~\ref{main} and the proofs for
auxiliary lemmas and the main result  will be given in
section~\ref{proof}.

\section{Main Result}\label{main}

Consider the $p_n\times p_n$ submatrix $\mathbf{A}$, truncated from a $n\times n$ circular unitary matrix $\mathbf{U}$
in section~\ref{intro}.
Denote the $p_n$ eigenvalues of $\mathbf{A}$  as $\bz_1,\cdots, \bz_{p_n}$ with the joint density function given by \eqref{Good_forever}.

The limiting distribution for  the spectral radius $\max_{1\le j\le
p_n}|\bz_j|$ has been obtained by Jiang and Qi~\cite{JiangQi2017}
and Gui and Qi~\cite{GQ2018} under each of the following conditions:

\begin{equation}\label{c0}
0<h_1<\frac{p_n}{n}<h_2<1, \mbox{ where $h_1$ and  $h_2$ are two constants;}
\end{equation}

\begin{equation}\label{c1}
p_n\to\infty \mbox{ and } \frac{p_n}{n}\to 0\mbox{ as }n\to\infty;
\end{equation}

\begin{equation}\label{c2}
\frac{n-p_n}{(\log n)^3}\to\infty\mbox{ and } \frac{n-p_n}{n}\to 0\mbox{ as }n\to\infty;
\end{equation}

\begin{equation}\label{c3}
n-p_n\to\infty \mbox{ and } \frac{n-p_n}{\log n}\to 0\mbox{ as }n\to\infty;
\end{equation}

\begin{equation}\label{c4}
n-p_n=k\ge 1 \mbox{ is a fixed integer.}
\end{equation}

Theorems~\ref{thm1} to \ref{thm3} below are summarized from Jiang
and Qi~\cite{JiangQi2017} and Gui and Qi~\cite{GQ2018}. The main
contribution of the present paper is Theorem~\ref{thm-main}.

\begin{theorem}\label{thm1}
Assume that $\bz_1, \cdots, \bz_{p_n}$ have density as in
(\ref{Good_forever}), and  $\{p_n\}$ is a sequence of positive
integers satisfying $1\le p_n< n$ and
\begin{equation}\label{combined}
 p_n\to\infty\mbox{ and } \frac{n-p_n}{(\log n)^3}\to\infty \mbox{ as }n\to\infty.
\end{equation}
Then $(\max_{1\leq j \leq p_n}|\bz_j|-A_n)/B_n$ converges weakly to
the Gumbel distribution $\Lambda(x)=\exp(-e^{-x})$, $x\in
\mathbb{R}$,
where
$A_n=c_n+\frac{1}{2}(1-c_n^2)^{1/2}(n-1)^{-1/2}a_n$,
$B_n=\frac{1}{2}(1-c_n^2)^{1/2}(n-1)^{-1/2}b_n$,
\[
c_n=\Big(\frac{p_n-1}{n-1}\Big)^{1/2},~~ b_n=b\Big(\frac{nc_n^2}{1-c_n^2}\Big), ~~a_n=a\Big(\frac{nc_n^2}{1-c_n^2}\Big)
\]
with \beaa a(y)=(\log y)^{1/2}-(\log y)^{-1/2}\log(\sqrt{2\pi}\log
y)\ \mbox{ and } \ b(y)=(\log y)^{-1/2} ~\mbox{ for }y>3. \eeaa
\end{theorem}

%%%where $k_n=n-p_n$.

%%%Gumbel-type, Fr\'echet-type, and

\begin{theorem}\label{thm2}
Under condition \eqref{c3},  $(\max_{1\leq j \leq p_n}|\bz_j|-A_n)/B_n$ converges weakly to the Gumbel distribution
$\Lambda(x)=\exp(-e^{-x})$, $x\in \mathbb{R}$, where $A_n=(1-a_n/n)^{1/2}$ and $B_n=a_n/(2nk_n)$ with $k_n=n-p_n$, and $a_n$ is given by
\[
\frac{1}{(k_n-1)!}\int^{a_n}_0t^{k_n-1}e^{-t}dt=\frac{k_n}{n}.
\]

\end{theorem}

\begin{theorem}\label{thm3}
Under condition \eqref{c4},   $\frac{2n^{1+1/k}}{((k+1)!)^{1/k}}(\max_{1\leq j \leq p_n}|\bz_j|-1)$ converges weakly to
the reversed Weibull distribution $W_{k}(x)$ defined as
\[
W_{k}(x)=\left\{
                \begin{array}{ll}
                 \exp(-(-x)^{k}) , & \hbox{$x\le 0$;} \\
                  1, & \hbox{$x>0$.}
                \end{array}
              \right.
\]
\end{theorem}

\noindent{\bf Remark 1.} Theorems~\ref{thm2} and \ref{thm3} are
proved in Gui and Qi~\cite{GQ2018}. Theorem~\ref{thm1} reduces to
Theorem 2 in Jiang and Qi~\cite{JiangQi2017} under \eqref{c0} and to
Theorem 2 in Gui and Qi~\cite{GQ2018} under \eqref{c1} or
\eqref{c2}. Note that condition \eqref{combined} combines conditions
\eqref{c0}, \eqref{c1} and \eqref{c2}.  In fact, Theorem~\ref{thm1}
can be concluded from Theorem 2 in Jiang and Qi~\cite{JiangQi2017}
and Theorem 2 in Gui and Qi~\cite{GQ2018} by using subsequence
arguments. A proof can be outlined as follows.   Let $\{p_n\}$ be
any sequence satisfying \eqref{combined}.  Then for any subsequence
$\{n'\}$ of positive integers, there always exists its further
subsequence, say $\{n''\}$, such that one of the three conditions
\eqref{c0}, \eqref{c1} and \eqref{c2} holds along the subsequence
$\{n''\}$. By applying Theorem 2 in Jiang and Qi~\cite{JiangQi2017}
or Theorem 2 in Gui and Qi~\cite{GQ2018}, we know that
Theorem~\ref{thm1} holds along the subsequence $\{n''\}$. This is
sufficient to conclude Theorem~\ref{thm1} above.

\vspace{10pt}

When $n-p_n$ is of order between $\log n$ and $(\log n)^3$, neither
of conditions from \eqref{c0} to \eqref{c4} holds. In this paper,
we consider the following condition
\begin{equation}\label{cc}
k_n=n-p_n\to\infty~\mbox{ and } \frac{k_n(\log n)^3}n\to 0~~\mbox{
as }~n\to\infty.
\end{equation}
The range of $p_n$ here is wide enough to cover the gap that is not considered in Jiang and Qi~\cite{JiangQi2017} and Gui and Qi~\cite{GQ2018}.

To define the normalizing constants for $\max_{1\leq j \leq p_n}|\bz_j|$,  set $\lambda_n$ as the solution to
\begin{equation}\label{gn-lambda}
g_n(\lambda):=\lambda-1-\log(\lambda)+\frac{2}{k_n}\log(1-\lambda) =\frac{1}{k_n}\log(\frac{n}{2\pi k_n^{3/2}})
\end{equation}
in $(0,1)$.  We see that $g_n(\lambda)$ is decreasing in $(0,1)$ by noting
\[
g_n'(\lambda)=1-\frac{1}{\lambda}-\frac{2k_n}{1-\lambda}<0 ~~~\mbox{ for } \lambda\in (0,1).
\]
Since $g_n(0+)=\infty$ and $g_n(1-)=-\infty$, a unique solution to $g_n(\lambda)=c$ in $(0,1)$ exists for any constant $c$.

Our main contribution in the paper is the following Theorem~\ref{thm-main}, which confirms the conjecture by Gui and Qi~\cite{GQ2018}.

\begin{theorem}\label{thm-main}
Under condition \eqref{cc},  $(\max_{1\leq j \leq p_n}|\bz_j|-A_n)/B_n$ converges weakly to the Gumbel distribution
$\Lambda(x)=\exp(-e^{-x})$, $x\in \mathbb{R}$, where
\[
A_n=(1-\frac{k_n\lambda_n}{n})^{1/2},~~~~B_n=\frac{\lambda_n}{2A_nn(1-\lambda_n)}.
\]
\end{theorem}

\vspace{10pt}

\noindent{\bf Remark 2.} The eigenvalues for truncation of the
circular unitary ensemble form a determinantal point process and
share the property of intrinsic independence. This property is very
helpful in investigating  both the asymptotic distribution of the
spectral radius and the empirical spectral distribution of the
eigenvalues from a determinantal point process; see, e.g., Jiang and
Qi~\cite{JiangQi2019}, Chang and Qi~\cite{ChangQi2017},  Chang, Li
and Qi~\cite{ChangLiQi2018} for more work on limiting empirical
spectral distributions for non-Hermitian random matrices. The proof
of Theorem~\ref{thm-main} is quite lengthy and will be split into a
series of auxiliary lemmas in section~\ref{proof}. For the case
under heavy truncation,  Gui and Qi~\cite{GQ2018} and Jiang and
Qi~\cite{JiangQi2017} employ moderate deviation principles for sum
of independent random variable, but this approach does not work any
more for our case. In fact, when $k_n=n-p_n$ is of order between
$\log n$ and $(\log n)^3$, we need a uniform estimate of the
probability for a Gamma($k_n$) random variable falling into the
interval $(0,x]$, where $x$ is between $0$ and some constant $c_n$
with $c_n<k_n$. Obviously, this is beyond the range one can apply
moderate deviation principles for Gamma($k_n$) since Gamma($k_n$) is
the sum of $k_n$ independent Gamma($1$) random variables. Instead,
we obtain a fine estimate in Lemma~\ref{Emme-lemma} below for
large-parameter incomplete Gamma function via a result in
Temme~\cite{Temme}. This lemma, together with
Lemmas~\ref{boundoflambda}, \ref{lambdaunif} and \ref{onetermilimit}
on several estimates of functions of the solution $\lambda_n$ to
equation \eqref{gn-lambda}, enables us to prove
Lemma~\ref{sumofprobability} and Theorem~\ref{thm-main}. Meanwhile,
this method may not be easily extended to prove the results from
Jiang and Qi~\cite{JiangQi2017} and Gui and Qi~\cite{GQ2018} in
general since approximations for some other terms will get worse if
$k_n$ is too large. Fortunately, the range of $p_n$ under condition
\eqref{cc} is wide enough to bridge the gap in the literature.

\section{Proofs}\label{proof}

We need the following notation in our proofs. We use the symbol
$C_n\sim D_n$ to denote the relationship
$\lim_{n\to\infty}\frac{C_n}{D_n}=1$. For random variables $\{X_n,
\, n\geq 1\}$ and constants $\{a_n, \, n\geq 1\}$, we write
$X_n=O_p(a_n)$ if $\lim_{x\to
+\infty}\limsup_{n\to\infty}P(|\frac{X_n}{a_n}|\geq x)=0$, and we
write $X_n=o_p(a_n)$ if $\frac{X_n}{a_n}\to 0$ in probability. It is
well known that $\frac{X_n}{a_nb_n}\to 0$ in probability as
$n\to\infty$ if $X_n=O_p(a_n)$ and $\{b_n, \, n\geq 1\}$ is a
sequence of constants with $\lim_{n\to\infty}b_n=\infty$.

As in Gui and Qi~\cite{GQ2018}, assume that $\{U_i, ~i\ge 1\}$ is a
sequence of independent and identically distributed (i.i.d.) random
variables uniformly distributed over $(0,1)$, and $U_{1:n}\le
U_{2:n}\le \cdots\le U_{n:n}$ denote the order statistics of $U_1,
U_2, \cdots, U_n$ for each $n\ge 1$. Then $U_{i:n}$ has a Beta($i,
n-i+1$) distribution with density function given by
\[
f_{i:n}(x)=\frac{n!}{(i-1)!(n-i)!}
x^{i-1}(1-x)^{n-i}, ~~~0<x<1,
\]
and its cumulative distribution function (cdf) is denoted by $F_{i:n}(x)$, $0\le x\le 1$.

For each $n\geq 2$, let $\{Y_{nj},\, 1\leq j \leq p_n\}$ be independent random variables such that
$Y_{nj}$ and $(U_{p_n+1-j: n-j})^{1/2}$ have the same distribution for each $j$.  Jiang and Qi~\cite{JiangQi2017}  have shown that
$\max_{1\leq j \leq p_n}|\bz_j|^2$ and $\max_{1\leq j \leq p_n}Y_{nj}^2$ have the same distribution.

Next, we express each Beta random variable in terms of Gamma random
variables. From equation (2.2.1) on page 12 in Ahsanullah and
Nevzorov~\cite{AN2015}, we have, for each $1\le k\le n$,  $U_{k:n}$
and $\sum^k_{i=1}E_i/\sum^{n+1}_{i=1}E_i$ have the same
distribution, where $\{E_i,~i\ge 1\}$ is a sequence of independent
random variables with the standard exponential distribution. In
fact,  if we assume that $\{E_{ij}, i\ge 1, j\ge 1\}$ are
independent random variables with the standard exponential
distribution, then
$\{\sum^{p_n+1-j}_{i=1}E_{i,j}/\sum^{n+1-j}_{i=1}E_{i,j},~1\le j\le
p_n\}$ are independent random variables, and for each $1\le j\le
p_n$, $\sum^{p_n+1-j}_{i=1}E_{i,j}/\sum^{n+1-j}_{i=1}E_{i,j}$ and
$U_{p_n+1-j: n-j}$ are identically distributed, which implies that
$\{\sum^{p_n+1-j}_{i=1}E_{i,j}/\sum^{n+1-j}_{i=1}E_{i,j},~1\le j\le
p_n\}$ and $\{Y_{nj}^2,~1\le j\le p_n\}$ are identically
distributed.  For simplicity, we assume
\begin{equation}\label{Ynj2}
Y_{nj}^2=\frac{\sum^{p_n+1-j}_{i=1}E_{i,j}}{\sum^{n+1-j}_{i=1}E_{i,j}}=1-\frac{S_j}{T_{n+1-j}},\mbox{
for }1\le j\le p_n
\end{equation}
where $S_j=\sum^{n+1-j}_{i=p_n-j+2}E_{i,j}$ and $T_{n+1-j}=\sum^{n+1-j}_{i=1}E_{i,j}$.   Then we have
\begin{equation}\label{identical}
P\Big(\max_{1\leq j \leq p_n}|\bz_j|^2\le t\Big)=P\Big(\max_{1\leq j \leq p_n}Y_{nj}^2\le t\Big)=\prod^{p_n}_{j=1}F_{p_n+1-j:n-j}(t)
\end{equation}
for $0<t<1$, and
\[
1-F_{1:k_n}(x)\le 1-F_{2: k_n+1}(x)\le \cdots\le 1-F_{p_n: n-1}(x)
\]
for $x\in (0,1)$.  See the proof of Theorem 2 in Jiang and Qi~\cite{JiangQi2017}.

It is easily seen that $\{S_j, ~1\le j\le p_n\}$ are i.i.d. random variables with Gamma ($k_n$) distribution.

%%Beta and F distribution:
%%Beyer, W. H. CRC Standard Mathematical Tables, 28th ed. Boca Raton, FL: CRC Press, p. 536, 1987.

%%We will need some properties of Gamma distributions.  Assume $\alpha>0$ and $\beta>0$. The density of a Gamma($\alpha, \beta$)
%%distribution is given by
%%\[
%%g(x, \alpha, \beta)=\frac{1}{\beta^\alpha\Gamma(\alpha)}x^{\alpha-1}e^{-x/\beta},~~~x>0,
%%\]
%%where $\Gamma(\alpha)$ is the gamma function defined as
%%\[
%%\Gamma(\alpha)=\int^\infty_0x^{\alpha-1}e^{-x}dx.
%%\]

%\subsection{Preliminary Lemmas}

We will present some useful lemmas before we prove our main
result.

\begin{lemma}\label{prod2sum} (Gui and Qi~\cite{GQ2018})
Suppose $\{l_n,\, n\geq 1\}$ is a sequence of positive integers. Let
$z_{nj}\in [0,1)$ be real numbers for $1\leq j \leq l_n$ such that
$\max_{1\le j\le l_n}z_{nj}\to 0$ as $n\to\infty$. Then
$\displaystyle\lim_{n\to\infty}\prod^{l_n}_{j=1}(1-z_{nj})\in (0,1)$
exists if and only if the limit
$\displaystyle\lim_{n\to\infty}\sum^{l_n}_{j=1}z_{nj}=:z\in
(0,\infty)$ exists and the relationship of the two limits is given
by
\[
\lim_{n\to\infty}\prod^{l_n}_{i=1}(1-z_{ni})=e^{-z}.
\]
\end{lemma}

%\hfill$\blacksquare$

%\begin{lemma}\label{little} (Gui and Qi, 2018) Let $\{l_n\}$ be a sequence of positive integers such that $l_n\to\infty$ and for each $n$, %$\{z_{nj}, ~1\le j\le l_n\}$ are non-negative numbers such that $z_{nj}$ is non-increasing in $j$ with $z_{n1}>0$.  Then for any sequence of %positive integers $\{r_n\}$satisfying that $r_n<l_n$ for all large $n$ and $r_n/l_n\to 1$ as $n\to\infty$, we have
%\[
%\frac{\sum^{l_n}_{j=1}z_{nj}}{\sum^{r_n}_{j=1}z_{nj}}\to 1
%\]
%as $n\to\infty$.
%\end{lemma}

%\hfill$\blacksquare$

\begin{lemma}\label{less} (Gui and Qi~\cite{GQ2018}) Assume that $1\le p_n<n$ and $p_n\to \infty$ as $n\to\infty$.
Let $\{r_n\}$ be a sequences of integers such that $r_n<p_n$ and
$p_n/r_n\to 1$ as $n\to \infty$.   Assume $\alpha_n>0$ and $\beta_n$
are real numbers such
that
$\lim_{n\to\infty}P(Y_{n1}^2>\beta_n+\alpha_nx)= 0$ for any
$x\in \mathbb{R}$. If $(\max_{1\le j\le
r_n}Y_{nj}^2-\beta_n)/\alpha_n$ converges in distribution to a cdf
$G$, then $(\max_{1\le j\le p_n}Y_{nj}^2-\beta_n)/\alpha_n$
converges in distribution to  the same distribution $G$.
\end{lemma}

\begin{lemma}\label{bird} (Gui and Qi~\cite{GQ2018})
Let $Z_n$ be nonnegative random variables such that $(Z_n^2-\beta_n)/\alpha_n$ converges weakly to a cdf $G(x)$, where
$\alpha_n>0$ and $\beta_n>0$ are constants satisfying that $\lim_{n\to\infty}\alpha_n/\beta_n=0$. Then
\[
\frac{Z_n-\beta_n^{1/2}}{\alpha_n/(2\beta_n^{1/2})} \mbox{ converges weakly to } G.
\]
\end{lemma}

%(Gui and Qi, 2018)
%Define
%\begin{equation}\label{vj}
%V_{p_n-j+1:n-j}=\frac{(n-j)^{3/2}}{((p_n-j)(n-p_n))^{1/2}}\Big(U_{p_n-j+1:n-j}-\frac{p_n-j}{n-j}\Big).
%\end{equation}
%Assume that $k_n=n-p_n\to\infty$ and $k_n/n\to 0$ as $n\to\infty$.    Then for any $\delta_n>0$
%such that $\delta_n\to\infty$ and $\delta_n=o(k_n^{1/6})$
%\begin{equation}\label{eq1}
%P(V_{p_n-j+1: n-j}>x)=(1+o(1))(1-\Phi(x))
%\end{equation}
%uniformly over $0\le x\le \delta_n$, $1\le j\le p_n-k_n$ as $n\to\infty$.

\vspace{10pt}

Recall we just define that $\{S_j,~ 1\le k\le p_n\}$ are i.i.d. Gamma($k_n$) random variables. For convenience, we assume that
 $\{S_j,~ j\ge 1\}$ are i.i.d. random variables with Gamma($k_n$)
 distribution.

Define the cumulative distribution function (cdf) of Gamma($a$) random variable (incomplete gamma function) as
\[
P(a,z)=\frac{1}{\Gamma(a)}\int^z_0x^{a-1}e^{-x}dx,~~~z\ge 0,
\]
and the error function
\[
\text{erfc}(z)=\frac{2}{\sqrt{\pi}}\int^\infty_ze^{-x^2}dx, ~~~z\ge 0.
\]

Write
\[
\tau(\lambda)=\lambda-1-\log\lambda, ~~~\lambda>0.
\]
It is easy to see that for $\lambda>0$
\[
\tau(\lambda)=\lambda-1-\log \lambda=(1-\lambda)^2\int^1_0\frac{s}{1-(1-\lambda)s}ds.
\]
We see that $\tau(\lambda)\ge 0$ for $\lambda>0$. Since $\min(\lambda,1)\le 1-(1-\lambda)s\le\max(1, \lambda)$ for $0<s<1$, we have
\begin{equation}\label{tau0}
\frac{1}{2}\frac{(1-\lambda)^2}{\max(1,\lambda)}\le \tau(\lambda)\le \frac{1}{2}\frac{(1-\lambda)^2}{\min(1,\lambda)} ~~~\lambda>0
\end{equation}
and conclude that
\begin{equation}\label{tau}
\sqrt{\tau(\lambda)}\ge \frac{1-\lambda}{\sqrt{2}}, ~~~0<\lambda<1.
\end{equation}

We can also verify that for $s, t>0$
\begin{equation}\label{taust}
\tau(st)=\tau(s)+\tau(t)+(s-1)(t-1).
\end{equation}
This property will be used later.

Define
\begin{equation}\label{phialambda}
\phi(a, \lambda)=\frac{1}{\sqrt{2\pi a}}e^{-a\tau(\lambda)}=\frac{1}{\sqrt{2\pi a}}e^{-a(\lambda-1-\log \lambda)}, ~~~~\lambda>0, ~a>0.
\end{equation}

\begin{lemma}\label{Emme-lemma} Let $\delta_n$ be a sequence of positive numbers
such that $\delta_n\to\infty$ and $\delta_n/\sqrt{k_n}\to 0$ as $n\to\infty$. Then
\begin{equation}\label{repPP}
P(k_n,k_n\lambda)=(1+o(1))\frac{1}{\sqrt{2\pi k_n}(1-\lambda)}\exp(-k_n\tau(\lambda))
\end{equation}
uniformly over $0<\lambda\le 1-\delta_n/\sqrt{k_n}$ as $n\to\infty$.
\end{lemma}

\noindent{\it Proof.}  It follows from  equations (2.15) and (4.3) in  Temme~\cite{Temme} that
\begin{equation}\label{temme-app}
|P(a, a\lambda)-\frac12\text{erfc}(-\eta\sqrt{\frac{a}{2}})+\frac{c_0(\lambda)}{\sqrt{2\pi a}}\exp(-\frac{1}{2}a\eta^2)|\le
\frac{C}{a\sqrt{2\pi a}}\exp(-\frac{1}{2}a\eta^2)+\frac{Ce^{a}a^{-a}\Gamma(a)}{a\sqrt{2\pi a}} P(a,a\lambda)
\end{equation}
holds uniformly for $0<\lambda<1$ and $a>0$,  where $C>0$ is a universal constant, $c_0(\lambda)=\frac{1}{\lambda-1}-\frac{1}{\eta}$,
and
\begin{equation}\label{etatau}
\eta=-(2(\lambda-1-\log \lambda))^{1/2}=-\sqrt{2\tau(\lambda)}
~~~\mbox{ for } 0<\lambda<1.
\end{equation}

%%%%-\eta=\sqrt{2\tau(\lambda)},~~~\frac{1}{2}\eta^2=\tau(\lambda).

From Stirling's formula, see, e.g., Formula 6.1.38 in Abramowitz and Stegun~\cite{Abramowitz1972}
\[
\Gamma(a+1)=\sqrt{2\pi}a^{a+\frac12}\exp(-a+\frac{\theta}{12a}),~~~\theta\in (0,1)
\]
we have
\[
e^{a}a^{-a}\Gamma(a)=e^{a}a^{-a-1}\Gamma(a+1)\le \sqrt{2\pi}e^{1/12}a^{-1/2}
\]
for all $a\ge 1$, which together with \eqref{temme-app} implies that for some universal constant $C>0$
\[
|P(a,a\lambda)-\frac12\text{erfc}(-\eta\sqrt{\frac{a}{2}})+\frac{c_0(\lambda)}{\sqrt{2\pi a}}\exp(-\frac{1}{2}a\eta^2)|\le
\frac{C}{a\sqrt{2\pi a}}\exp(-\frac{1}{2}a\eta^2)+\frac{C}{a^2}P(a,a\lambda)
\]
holds uniformly for $0<\lambda<1$ and $a>1$.  By setting
\[
c(a, \lambda)=\frac{P(a,a\lambda)-\frac12\text{erfc}(-\eta\sqrt{\frac{a}{2}})+\frac{c_0(\lambda)}{\sqrt{2\pi a}}\exp(-\frac{1}{2}a\eta^2)}{
\frac{1}{a\sqrt{2\pi a}}\exp(-\frac{1}{2}a\eta^2)+\frac{1}{a^2}P(a,a\lambda)},
\]
we have  $|c(a, \lambda)|\le C$ for $0<\lambda<1$ and $a>1$, and
\begin{equation}\label{repP}
P(a,a\lambda)-\frac12\text{erfc}(-\eta\sqrt{\frac{a}{2}})+\frac{c_0(\lambda)}{\sqrt{2\pi a}}\exp(-\frac{1}{2}a\eta^2)=c(a,\lambda)(
\frac{1}{a\sqrt{2\pi a}}\exp(-\frac{1}{2}a\eta^2)+\frac{1}{a^2}P(a,a\lambda)).
\end{equation}

Assume $x>0$. Since
\[
\int^{\infty}_x\frac{e^{-t^2}}{t^2}dt<\frac1{x^3}\int^{\infty}_xte^{-t^2}dt=\frac{1}{2x^3}\int^\infty_xe^{-t^2}dt^2
=\frac{e^{-x^2}}{2x^3}
\]
and by using integration by parts
\begin{eqnarray*}
\int^\infty_xe^{-t^2}dt&=&-\frac12\int^{\infty}_x\frac{de^{-t^2}}{t}\\
&=&\frac{1}{2}\frac{e^{-x^2}}{x}-\frac12\int^{\infty}_x\frac{e^{-t^2}}{t^2}dt,
\end{eqnarray*}
we get
\begin{eqnarray*}
   (\frac{1}{2x}-\frac{1}{4x^3})e^{-x^2}<\int^\infty_xe^{-t^2}dt<\frac{1}{2x}e^{-x^2}
\end{eqnarray*}
%i.e.,
%\begin{eqnarray*}
%   (\frac{1}{\sqrt{2}x}-\frac{1}{\sqrt{2}x^3})e^{-x^2/2}<\int^\infty_{x/\sqrt{2}}e^{-t^2}dt<\frac{1}{\sqrt{2}x}e^{-x^2/2}
%\end{eqnarray*}
Therefore, we can define function $h(x)$ such that
\begin{equation}\label{erfc}
 \frac{1}{\sqrt{\pi}}\int^\infty_{x}e^{-t^2}dt=(\frac{1}{\sqrt{2}x}
 -\frac{h(x)}{\sqrt{2}x})\frac{1}{\sqrt{2\pi}}e^{-x^2},
\end{equation}
where $0<h(x)<\frac{1}{2x^{2}}$ for $x>0$.

By using \eqref{erfc} and \eqref{etatau}, we have
\begin{eqnarray*}
&&\frac12\text{erfc}(-\eta\sqrt{\frac{a}{2}})-\frac{c_0(\lambda)}{\sqrt{2\pi a}}\exp(-\frac{1}{2}a\eta^2)\\
&=&\frac{1}{2}\text{erfc}(\sqrt{a\tau(\lambda)})-(\frac{1}{\sqrt{2\tau(\lambda)}}-\frac1{1-\lambda})\frac{1}{\sqrt{2\pi a}}\exp(-a\tau(\lambda))\\
&=&\frac1{1-\lambda}\Big(1-\frac{1-\lambda}{\sqrt{2\tau(\lambda)}}h(\sqrt{a\tau(\lambda)})\Big)\frac{1}{\sqrt{2\pi a}}\exp(-a\tau(\lambda))\\
&=&\Big(1-\frac{1-\lambda}{\sqrt{2\tau(\lambda)}}h(\sqrt{a\tau(\lambda)})\Big)\frac{\phi(a,\lambda)}{1-\lambda}.
\end{eqnarray*}
Then it follows from \eqref{repP} that
\[
(1-\frac{c(a,\lambda)}{a^2})P(a,a\lambda)=\Big(1-\frac{1-\lambda}{\sqrt{2\tau(\lambda)}}h(\sqrt{a\tau(\lambda)})+\frac{c(a,\lambda)(1-\lambda)}{a}\Big)\frac{\phi(a,\lambda)}{1-\lambda}
\]
and thus
\[
P(a,a\lambda)=\frac{\Big(1-\frac{1-\lambda}{\sqrt{2\tau(\lambda)}}h(\sqrt{a\tau})+\frac{c(a,\lambda)(1-\lambda)}{a}\Big)}
{1-\frac{c(a,\lambda)}{a^2}}\frac{\phi(a,\lambda)}{1-\lambda}
\]
uniformly over $0<\lambda<1$ and $a>1$.

Now let $a=k_n$. Let $\delta_n$ be any sequence of positive numbers such that $\delta_n\to\infty$ and $\delta_n/\sqrt{k_n}\to 0$ as $n\to\infty$.  From
\eqref{tau} we have
\[
\sqrt{k_n\tau(\lambda)}\ge \sqrt{\frac{k_n}2}(1-\lambda)\ge \frac{\delta_n}{\sqrt{2}}\to\infty
\]
if $0<\lambda\le 1-\delta_n/\sqrt{k_n}$, which implies $h(\sqrt{k_n\tau(\lambda)})\to 0$ uniformly over $0<\lambda\le 1-\delta_n/\sqrt{k_n}$ as $n\to\infty$.  Therefore, we conclude that
%\begin{equation}\label{repPP}
%P(k_n,k_n\lambda)=(1+o(1))\frac{\phi(a,\lambda)}{1-\lambda}=(1+o(1))\frac{1}{\sqrt{2\pi k_n}(1-\lambda)}\exp(-(\lambda-1-\log\lambda))
%\end{equation}
\[
P(k_n,k_n\lambda)=(1+o(1))\frac{\phi(a,\lambda)}{1-\lambda}=(1+o(1))\frac{1}{\sqrt{2\pi k_n}(1-\lambda)}\exp(-k_n\tau(\lambda))
\]
uniformly over $0<\lambda\le 1-\delta_n/\sqrt{k_n}$ as $n\to\infty$,
i.e. \eqref{repPP} holds.  This completes the proof of the lemma.
\hfill $\blacksquare$

\begin{lemma}\label{boundoflambda} Under condition \eqref{cc} we have
\begin{equation}\label{1minuslambda}
\sqrt{k_n}(1-\lambda_n)\to\infty~~~\mbox{ as }n\to\infty
\end{equation}
and
\begin{equation}\label{1-lambda-upper}
\sqrt{k_n}(1-\lambda_n)=O(\sqrt{\log n})~~~\mbox{ as }n\to\infty.
\end{equation}
\end{lemma}

\noindent{\it Proof.} Since
\[
g_n(\lambda)=\lambda-1-\log(\lambda)+\frac{2}{k_n}\log(1-\lambda)
\]
is decreasing in $\lambda\in (0,1)$, we have for any $\delta>0$,
$\sqrt{k_n}(1-\lambda_n)>\delta$ if and only if
$g_n(1-\delta/\sqrt{k_n})<\frac{1}{k_n}\log(\frac{n}{2\pi
k_n^{3/2}})$.
To prove \eqref{1minuslambda},  it suffice to show that
for any $\delta>0$,
$g_n(1-\delta/\sqrt{k_n})<\frac{1}{k_n}\log(\frac{n}{2\pi
k_n^{3/2}})$ for all large $n$. In fact, for any fixed $\delta>0$,
we have from \eqref{tau0} that for all large $n$
\begin{eqnarray*}
g_n(1-\frac{\delta}{\sqrt{k_n}})&=&\tau(1-\frac{\delta}{\sqrt{k_n}})+\frac{2}{k_n}\log(\frac{\delta}{\sqrt{k_n}})\\
&\le& \frac{1}{2}\frac{\delta^2}{k_n(1-\delta/\sqrt{k_n})}+\frac{2}{k_n}\log(\frac{\delta}{\sqrt{k_n}})\\
&\le& \frac{\delta^2}{k_n}+\frac{2}{k_n}\log(\frac{\delta}{\sqrt{k_n}})\\
&<& \frac{1}{k_n}\log(\frac{n}{2\pi \delta^2\sqrt{k_n}})+\frac{1}{k_n}\log(\frac{\delta^2}{k_n})\\
&=&\frac{1}{k_n}\log(\frac{n}{2\pi k_n^{3/2}}),
\end{eqnarray*}
proving \eqref{1minuslambda}.

Now we prove \eqref{1-lambda-upper}. By using \eqref{tau}, we have
\[
(1-\lambda_n)^2\le
2\tau(\lambda_n)=2g(\lambda_n)-\frac4{k_n}\log(1-\lambda_n)=\frac{2}{k_n}\log(\frac{n\sqrt{k_n}}{2\pi
})-\frac{4}{k_n}\log(k_n(1-\lambda_n)).
\] From \eqref{1minuslambda},
we have $\log(k_n(1-\lambda_n))>0$ for all large $n$. Therefore, we
get
\[
k_n(1-\lambda_n)^2<2\log(\frac{n\sqrt{k_n}}{2\pi})=O(\log n)~~~\mbox{ as }n\to\infty,
\]
which proves \eqref{1-lambda-upper}. \hfill $\blacksquare$

\vspace{10pt}

For convenience, we will introduce more notations for the rest of the paper.

Define for $x\in \mathbb{R}$
\begin{equation}\label{lambdanx}
\lambda_n(x)=\lambda_n(1+\frac{x}{k_n(1-\lambda_n)}),
\end{equation}
and for $1\le j<n$
\begin{equation}\label{lambdanj}
\lambda_{n,j}(x)=\frac{n+1-j}{n}\lambda_n(x).
\end{equation}

\begin{lemma}\label{lambdaunif} Assume condition~\eqref{cc} holds.
We have for any fixed $x\in\mathbb R$ that
\begin{equation}\label{1-lambda2}
\frac{1-\lambda_n(x)}{1-\lambda_n}-1\to 0.
\end{equation}
If further we assume that $\{j_n\}$ is a sequence of positive integers with $1< j_n< n-1$ such that
\begin{equation}\label{knjn}
\frac{k_nj_n(1-\lambda_n)}n\to\infty ~\mbox{ and }~\frac{k_nj_n^2}{n^2}\to 0,
\end{equation}
then for any fixed $x\in\mathbb R$
\begin{equation}\label{1-lambda1}
\max_{1\le j\le j_n}|\frac{1-\lambda_{n,j}(x)}{1-\lambda_n}-1|\to 0.
\end{equation}
\end{lemma}

\noindent{\it Proof.} The proofs are omitted here since they are straightforward by using Lemma~\ref{boundoflambda} and given conditions.
\hfill $\blacksquare$

\begin{lemma}\label{onetermilimit} Assume condition \eqref{cc} holds. Then  with $\phi(a,\lambda)$ defined in \eqref{phialambda}
we have for $x\in \mathbb{R}$
\begin{equation}\label{first-term}
\frac{1}{(1-\lambda_n)^2}\phi(k_n, \lambda_n(x))=(1+o(1))\frac{k_ne^{x}}{n}~\mbox{ as }n\to\infty.
\end{equation}
\end{lemma}

\noindent{\it Proof.} Fix $x\in \mathbb{R}$.
It follows from \eqref{1minuslambda} and \eqref{1-lambda2} that
\[
\delta_n(x):=\sqrt{k_n}(1-\lambda_n(x))=\sqrt{k_n}(1-\lambda_n)(1+o(1))\to\infty.
\]
Then $0<\lambda_n(x)<1$ for all large $n$, which will be used in the proof of Lemma~\ref{sumofprobability}.
By using \eqref{taust} we have from \eqref{tau0}
\begin{eqnarray*}
\tau(\lambda_n(x))&=&\tau(\lambda_n)+\tau(1+\frac{x}{k_n(1-\lambda_n)})+(\lambda_n-1)\frac{x}{k_n(1-\lambda_n)}\\
&=&\tau(\lambda_n)+\tau(1+\frac{x}{k_n(1-\lambda_n)})-\frac{x}{k_n}\\
&=&g_n(\lambda_n)-\frac{2}{k_n}\log(1-\lambda_n)+\frac{x}{k_n}+\tau(1+\frac{x}{k_n(1-\lambda_n)})\\
&=&\frac{1}{k_n}\log(\frac{n}{\sqrt{2\pi} k_n^{3/2}})-\frac{2}{k_n}\log(1-\lambda_n)-\frac{x}{k_n}+O(\frac{1}{k_n^2(1-\lambda_n)^2}).
\end{eqnarray*}
Therefore, we obtain
\begin{eqnarray*}
&&\frac{1}{(1-\lambda_n)^2}\phi(k_n, \lambda_n(x))\\
&=&\frac{1}{\sqrt{2\pi k_n}(1-\lambda_n)^2}e^{-k_n\tau(\lambda_n(x))}\\
&=&\frac{1}{\sqrt{2\pi k_n}(1-\lambda_n)^2}\exp\Big\{-\log(\frac{n}{\sqrt{2\pi} k_n^{3/2}})+2\log(1-\lambda_n)+x+O(\frac{1}{k_n(1-\lambda_n)^2})\Big\}\\
&=&\frac{k_ne^{x}}{n}\exp\Big\{O(\frac{1}{k_n(1-\lambda_n)^2})\Big\}\\
&=&\frac{k_ne^{x}}{n}(1+o(1)).
\end{eqnarray*}
This completes the proof. \hfill $\blacksquare$

\begin{lemma}\label{sumofprobability}Fix $x\in\mathbb{R}$.  Assume $\{j_n\}$ is a sequence of positive integers with $1< j_n< n-1$ such that
\eqref{knjn} holds.
Then
\begin{equation}\label{limit1}
\sum^{j_n}_{j=1}P(k_n, \frac{n+1-j}{n}\lambda_n(x))\to e^{x}.
\end{equation}
Furthermore, if $\{q_n\}$ is a sequence of positive integers such that $1<j_n\le q_n< n-1$ , then
\begin{equation}\label{limit2}
\sum^{q_n}_{j=1}P(k_n, \frac{n+1-j}{n}\lambda_n(x))\to e^{x}.
\end{equation}
\end{lemma}

\noindent{\it Proof.}   It follows from \eqref{1minuslambda} and
\eqref{1-lambda2} that
$\delta_{n}(x)=\sqrt{k_n}(1-\lambda_n(x))\to\infty$. Since
$0<\lambda_{n,j}\le \lambda_n(x)\le 1-\delta_n/\sqrt{k_n}$ for all
$1\le j<n$, we have from \eqref{repPP} that

\begin{equation}\label{allj}
P(k_n,k_n\lambda_{n,j}(x))=(1+o(1))\frac{1}{\sqrt{2\pi k_n}(1-\lambda_{n,j}(x))}\exp(-k_n\tau(\lambda_{n,j}(x)))
\end{equation}
uniformly over $1\le j<n$ as $n\to\infty$.

Since $\lambda_{n,j}(x)=\frac{n+1-j}{n}\lambda_n(x)$, we have from \eqref{taust} that
\[
\tau(\lambda_{n,j}(x))=\tau(\lambda_n(x))+\tau(1-\frac{j-1}{n})+(1-\lambda_n(x))\frac{j-1}{n}.
\]
From \eqref{allj} we have
\begin{equation}\label{allj1}
P(k_n,k_n\lambda_{n,j}(x))=(1+o(1))\frac{\phi(k_n,
\lambda_n(x))}{1-\lambda_{n,j}(x)}\exp\{-k_n\tau(1-\frac{j-1}{n})-\frac{(j-1)(1-\lambda_n(x))k_n}{n}\}
\end{equation}
uniformly over $1\le j<n$ as $n\to\infty$.

Note that \eqref{1-lambda1} and \eqref{1-lambda2} hold under \eqref{knjn}. Then we have from \eqref{allj1} and \eqref{tau0}
that
\begin{eqnarray*}
P(k_n,k_n\lambda_{n,j}(x))&=&(1+o(1))\frac{\phi(k_n, \lambda_n(x))}{1-\lambda_n(x)}\exp\{-k_nO(\frac{j_n^2}{n^2})-\frac{(j-1)(1-\lambda_n(x))k_n}{n}\}\nonumber\\
&=&(1+o(1))\frac{\phi(k_n, \lambda_n(x))}{1-\lambda_n(x)}\exp\{o(1)-\frac{(j-1)(1-\lambda_n(x))k_n}{n}\}\nonumber\\
&=&(1+o(1))\frac{\phi(k_n, \lambda_n(x))}{1-\lambda_n(x)}\exp\{-\frac{(j-1)(1-\lambda_n(x))k_n}{n}\}
\end{eqnarray*}
uniformly over $1\le j<j_n$ as $n\to\infty$, which yields that
\begin{eqnarray*}
\sum^{j_n}_{j=1}P(k_n,k_n\lambda_{n,j}(x))&=&(1+o(1))\frac{\phi(k_n, \lambda_n(x))}{1-\lambda_n(x)}\sum^{j_n}_{j=1}\Big(\exp\{-\frac{(1-\lambda_n(x))k_n}{n}\}\Big)^{j-1}\\
&=&(1+o(1))\frac{\phi(k_n, \lambda_n(x))}{1-\lambda_n(x)}\frac{1-\exp\{-\frac{(1-\lambda_n(x))(j_n+1)k_n}{n}\}}{1-\exp\{-\frac{(1-\lambda_n(x))k_n}{n}\}}\\
&=&(1+o(1))\frac{\phi(k_n, \lambda_n(x))}{1-\lambda_n(x)}\frac{1-\exp\{-(1+o(1))\frac{(1-\lambda_n)(j_n+1)k_n}{n}\}}{\frac{(1-\lambda_n(x))k_n}{n}}\\
&=&(1+o(1))\frac{\phi(k_n, \lambda_n(x))}{(1-\lambda_n(x))^2}\frac{n}{k_n}\\
&=&(1+o(1))\frac{\phi(k_n, \lambda_n(x))}{(1-\lambda_n)^2}\frac{n}{k_n}\\
&=&(1+o(1))e^{x}
\end{eqnarray*}
from \eqref{first-term}.

Next,  we prove \eqref{limit2} when $q_n>j_n$.  Note  that $\tau(1-\frac{j-1}{n})\ge 0$, and $1-\lambda_{n,j}(x)\ge 1-\lambda_n(x)$ for $1\le j\le q_n$. Then from \eqref{allj1} we have
 \begin{eqnarray*}
\sum^{q_n}_{j=1}P(k_n,k_n\lambda_{n,j}(x))&\le &
(1+o(1))\frac{\phi(k_n, \lambda_n(x))}{1-\lambda_{n,j}(x)}\sum^{q_n}_{j=1}\exp\{--\frac{(j-1)(1-\lambda_n(x))k_n}{n}\}\\
&=&(1+o(1))\frac{\phi(k_n, \lambda_n(x))}{1-\lambda_{n,j}(x)}\sum^{q_n}_{j=1}\Big(\exp\{-\frac{(1-\lambda_n(x))k_n}{n}\}\Big)^{j-1}\\
&\le &(1+o(1))\frac{\phi(k_n, \lambda_n(x))}{1-\lambda_n(x)}\frac{1}{1-\exp\{-\frac{(1-\lambda_n(x))k_n}{n}\}}\\
&=&(1+o(1))\frac{\phi(k_n, \lambda_n(x))}{1-\lambda_n(x)}\frac{1}{\frac{(1-\lambda_n(x))k_n}{n}}\\
&=&(1+o(1))\frac{\phi(k_n, \lambda_n(x))}{(1-\lambda_n(x))^2}\frac{n}{k_n}\\
&=&(1+o(1))\frac{\phi(k_n, \lambda_n(x))}{(1-\lambda_n)^2}\frac{n}{k_n}\\
&=&(1+o(1))e^{x}
\end{eqnarray*}
from \eqref{first-term}, which together with \eqref{limit1}, implies
\eqref{limit2} since
\[
\sum^{j_n}_{j=1}P(k_n,k_n\lambda_{n,j}(x))\le \sum^{q_n}_{j=1}P(k_n,k_n\lambda_{n,j}(x)).
\]
This completes the proof of the lemma.
\hfill $\blacksquare$

\begin{lemma}\label{minimum}
Assume $x_1, \cdots, x_n$ are positive numbers for $n\ge 1$ and $\varepsilon_{n, i}$, $1\le i\le n$ are such that
$\varepsilon_n=\max_{1\le i\le n}|\varepsilon_{n, i}|<1$. Then
\[
  |\min_{1\le i\le n}x_i(1+\varepsilon_{n, i})-\min_{1\le i\le n}x_i|\le \varepsilon_n\min_{1\le i\le n}x_i.
\]
\end{lemma}

\noindent{\it Proof.} For $1\le i\le n$ we have
\[
   x_i(1-\varepsilon_n)\le x_i(1+\varepsilon_{n, i})\le x_i(1+\varepsilon_n),
   \]
and thus we obtain that
\[
  \min_{1\le i\le n} x_i(1-\varepsilon_n)\le \min_{1\le i\le n} x_i(1+\varepsilon_{n, i})\le \min_{1\le i\le n}x_i (1+\varepsilon_n),
   \]
which implies
\[
-\varepsilon_n\min_{1\le i\le n}x_i\le \min_{1\le i\le n}x_i(1+\varepsilon_{n, i})-\min_{1\le i\le n}x_i \le \varepsilon_n\min_{1\le i\le n}x_i.
\]
This completes the proof of the lemma.  \hfill $\blacksquare$

\begin{lemma}\label{expansion1} Under condition \eqref{cc} we have
\[
\lim_{n\to\infty}P(Y_{n1}^2>1-\frac{k_n}{n}\lambda_n(-x))=
\lim_{n\to\infty}P(Y_{n1}^2>1-\frac{k_n\lambda_n}{n}+\frac{\lambda_nx}{n(1-\lambda_n)})=0
\]
for each $x\in\mathbb{R}.$
\end{lemma}

\noindent{\it Proof.} By using expression \eqref{Ynj2} we have
$Y_{n1}^2=1-\frac{S_1}{T_n}$, where $S_1$ has a Gamma ($k_n$)
distribution and $T_n$ has a Gamma($n$) distribution. From the
central limit theorem we have
\[
V_{n1}:=\frac{S_1-k_n}{k_n^{1/2}}\xrightarrow{d} N(0,1),
~~~V_{n2}:=\frac{T_n-n}{n^{1/2}}\xrightarrow{d} N(0,1).
\]
Then we have
\[Y_{n1}^2=1-\frac{k_n+k_n^{1/2}V_{n1}}{n+n^{1/2}V_{n2}}=1-\frac{k_n}{n}\frac{1+V_{n1}/k_n^{1/2}}{1+V_{n2}/n^{1/2}}
=1-\frac{k_n}{n}(1+\frac{V_{n1}}{k_n^{1/2}}+O_p(n^{-1/2})),
\]
and thus we get
\[
V_{n3}:=\frac{n}{k_n^{1/2}}(Y_{n1}^2-1+\frac{k_n}{n})=V_{n1}+O_p(\sqrt{\frac{k_n}n})\xrightarrow{d}
N(0,1),
\]
 which yields
\[
P(Y_{n1}^2>1-\frac{k_n\lambda_n}{n}+\frac{\lambda_nx}{n(1-\lambda_n)})
=P(V_{n3}>\sqrt{k_n}(1-\lambda_n)+\frac{\lambda_nx}{\sqrt{k_n}(1-\lambda_n)})\to
0
\] since $\sqrt{k_n}(1-\lambda_n)\to \infty$ as $n\to\infty$ from
\eqref{1minuslambda}.  This completes of the proof. \hfill
$\blacksquare$

\begin{lemma}\label{max}
Let $\{T_j, ~j\ge 1\}$ be a sequence of random variables, and for each $j\ge 1$, $T_j$ has a Gamma ($j$) distribution with density function
$t^{j-1}e^{-t}I(t>0)/(j-1)!$.  Then
\[
\max_{m_n\le j\le n}|\frac{T_{j}}{j}-1|=O(\frac{\sqrt{\log n}}{\sqrt{m_n}}) ~~\mbox{ almost surely (a.s.)},
\]
where $m_n$ is any sequence of integers such that $1\le m_n<n$ and $m_n/(\log n)^3\to \infty$ as $n\to\infty$.
\end{lemma}

\noindent{\it Proof.} Set $\tau_n=\frac{\sqrt{\log n}}{\sqrt{m_n}}$. By Theorem 1 on page 217 in Petrov~\cite{Petrov},
\[
P(T_j-j>x\sqrt{j})=(1+o(1))(1-\Phi(x)), ~~~P(T_j-j<-x\sqrt{j})=(1+o(1))(1-\Phi(x))
\]
uniformly over $0\le x\le d_j$ as $j\to\infty$, where $\Phi$ is the
cumulative distribution function for the standard normal random
variable,  $d_j$ is any sequence of positive numbers with
$d_j=o(j^{1/6})$. By setting $x=4\sqrt{\log n}$ when $m_n\le j\le n$
and using the approximation
$1-\Phi(x)\sim\frac{1}{\sqrt{2\pi}x}e^{-x^2/2}=\frac{1}{\sqrt{2\pi}4n^8\sqrt{\log
n}}$ we conclude that
\begin{eqnarray*}
&&P(\max_{m_n\le j\le n}|\frac{T_{j}}{j}-1|>2\tau_n)\\
&\le& \sum^n_{j=m_n}P(|\frac{T_{j}}{j}-1|>2\tau_n)\\
&\le& \sum^n_{j=m_n}P(|T_{j}-j|>2j\tau_n)\\
&\le& \sum^n_{j=m_n}P(|T_j-j|>2\sqrt{\log n}\sqrt{j})\\
&\le& \sum^n_{j=m_n}P(T_j-j>2\sqrt{\log n}\sqrt{j})+\sum^n_{j=m_n}P(T_j-j<-2\sqrt{\log n}\sqrt{j})\\
&\sim& \frac{2(n-m_n+1)}{\sqrt{2\pi}4n^8\sqrt{\log n}},
\end{eqnarray*}
which implies that $\sum_{n>5}P(\max_{m_n\le j\le n}|\frac{T_{j}}{j}-1|>2\tau_n)<\infty$.  The lemma is proved by using the Borel-Cantelli lemma.
\hfill $\blacksquare$

\vspace{10pt}

From now on, we define $m_n=[k_n(\log n)^3]$, the integer part of $k_n(\log n)^3$. Then $m_n>k_n$ for all large $n$ and $m_n/(\log n)^3\to\infty$
as $n\to\infty$.

\begin{lemma}\label{maxSj}   Under condition \eqref{cc} we have
\begin{equation}\label{Ln}
L_n:=({\frac{\lambda_n}{n(1-\lambda_n)}})^{-1} (\min_{1\le j\le n-m_n}\frac{S_j}{n+1-j}-\frac{k_n\lambda_n}{n})\xrightarrow{d} \Lambda_1
\end{equation}
where $\Lambda_1(x)=1-\Lambda(-x)$, $x\in \mathbb{R}$.
\end{lemma}

\noindent{\it Proof.}   Fix $x\in \mathbb{R}$.  Recall $\lambda_n(x)$ and $\lambda_{n,j}(x)$ are defined in \eqref{lambdanx} and
\eqref{lambdanj}, respectively.  We have
\begin{eqnarray*}
&&P\Big(({\frac{\lambda_n}{n(1-\lambda_n)}})^{-1} (\min_{1\le j\le n-m_n}\frac{S_j}{n+1-j}-\frac{k_n\lambda_n}{n})\le x\Big)\\
&=&P\Big(\min_{1\le j\le n-m_n}\frac{S_j}{n+1-j}\le \frac{k_n\lambda_n}{n}+ \frac{\lambda_n}{n(1-\lambda_n)} x\Big)\\
&=&P\Big(\min_{1\le j\le n-m_n}\frac{S_j}{n+1-j}\le \frac{k_n}{n}\lambda_n(x)\Big)\\
&=&1-P\Big(\min_{1\le j\le n-m_n}\frac{S_j}{n+1-j}>\frac{k_n}{n}\lambda_n(x)\Big)\\
&=&1-\prod^{n-m_n}_{j=1}P\Big(\frac{S_j}{n+1-j}>\frac{k_n}{n}\lambda_n(x)\Big)\\
&=&1-\prod^{n-m_n}_{j=1}P\Big(S_j>k_n\lambda_{n,j}(x)\Big)\\
&=&1-\prod^{n-m_n}_{j=1}(1-z_{nj}),
\end{eqnarray*}
where $z_{nj}=P(k_n, k_n\lambda_{n,j}(x))$ for $1\le j\le n-m_n$. Note that $z_{nj}$ is non-increasing in $j$, and $\lambda_{n,1}(x)=\lambda_n(x)$. It follows from \eqref{allj1} with $j=1$, \eqref{first-term} and \eqref{1-lambda2} that
\[
z_{n1}=(1+o(1))\frac{\phi(k_n, \lambda_n(x))}{1-\lambda_n(x)}\le 1+o(1))\frac{\phi(k_n, \lambda_n(x))}{(1-\lambda_n(x))^2}= (1+o(1))\frac{k_ne^{x}}{n}\to 0
\]
as $n\to\infty$.

To apply Lemma~\ref{sumofprobability}, we define $\delta_n=k_n^{1/2}(1-\lambda_n)$.  Then $\delta_n\to\infty$ from \eqref{1minuslambda}.
Set $j_n=[n/\sqrt{k_n\delta_n}]$, the integer part of $n/\sqrt{k_n\delta_n}$. Then as $n\to\infty$
\[
\frac{k_nj_n(1-\lambda_n)}{n}=\frac{k_n^{1/2}j_n\delta_n}{n}\sim \sqrt{\delta_n}\to\infty
\]
and
\[
\frac{k_nj_n^2}{n^2}\sim \frac{1}{\delta_n}\to 0,
\]
i.e.  \eqref{knjn} holds. Obviously we have $n-m_n>j_n$ for all large $n$. Therefore,  by applying Lemma~\ref{sumofprobability} with
$q_n=n-m_n$, we have
\[
\sum^{n-m_n}_{j=1}z_{nj}\to e^{x},
\]
which coupled with Lemma~\ref{prod2sum} yields that $\prod^{n-m_n}_{j=1}(1-z_{nj})\to \exp(-e^{x})=\Lambda(-x)$ as $n\to\infty$.
Hence, we get
\[ P\Big(({\frac{\lambda_n}{n(1-\lambda_n)}})^{-1} (\min_{1\le j\le n-m_n}\frac{S_j}{n+1-j}-\frac{k_n\lambda_n}{n})\le x\Big)
\to 1-\Lambda(-x)=\Lambda_1(x),
\]
which proves the lemma. \hfill $\blacksquare$

\begin{lemma}\label{maxY2} Under condition \eqref{cc} we have
\[
(\frac{\lambda_n}{n(1-\lambda_n)})^{-1}(\max_{1\le j\le n-m_n}Y_{nj}^2-(1-\frac{k_n\lambda_n}{n}))\xrightarrow{d} \Lambda.
\]
\end{lemma}

\noindent{\it Proof.} From Lemma~\ref{max} we have  $\varepsilon_n:=\max_{m_n<j\le n}|\frac{T_j}{j}-1|=\max_{1\le j\le n-m_n}|\frac{T_{n+1-j}}{n+1-j}-1|=O(\frac{\sqrt{\log n}}{\sqrt{m_n}})\to 0$ a.s. as $n\to \infty$.
Then we have
\[
1-\varepsilon_n\le \frac{T_{n+1-j}}{n+1-j}\le
1+\varepsilon_n~~~\mbox{ uniformly for } 1\le j\le n-m_m
\]
for all large $n$, i.e.
\[
1-\frac{\varepsilon_n}{1-\varepsilon_n}\le
\frac1{1+\varepsilon_n}\le \frac{n+1-j}{T_{n+1-j}}\le
\frac{1}{1-\varepsilon_n}=1+\frac{\varepsilon_n}{1-\varepsilon_n}~~~\mbox{
uniformly for } 1\le j\le n-m_m
\]
for large $n$.  By writing $\varepsilon_{nj}= \frac{n+1-j}{T_{n+1-j}}-1$, we  have from \eqref{Ynj2} that
$Y_{nj}^2=1-\frac{S_{j}}{n+1-j}(1+\varepsilon_{nj})$, and thus
\[
\max_{1\le j\le n-m_n}Y_{nj}^2=1-\min_{1\le j\le n-m_n}\frac{S_{j}}{n+1-j}(1+\varepsilon_{nj}).
\]
Recall $L_n$ is defined in equation \eqref{Ln}.  The above equation, together with Lemma~\ref{minimum}, yields that for all large $n$
\begin{eqnarray}\label{Deltan}
\Delta_n:&=&|(\frac{\lambda_n}{n(1-\lambda_n)})^{-1}\Big(\max_{1\le j\le n-m_n}Y_{nj}^2-(1-\frac{k_n\lambda_n}{n})\Big)
+L_n|\\
&=&(\frac{\lambda_n}{n(1-\lambda_n)})^{-1}|\min_{1\le j\le n-m_n}\frac{S_{j}}{n+1-j}(1+\varepsilon_{nj})-\min_{1\le j\le n-m_n}\frac{S_{j}}{n+1-j}|\nonumber\\
&\le& \frac{\varepsilon_n}{1-\varepsilon_n}(\frac{\lambda_n}{n(1-\lambda_n)})^{-1}\min_{1\le j\le n-m_n}\frac{S_{j}}{n+1-j}. \nonumber
\end{eqnarray}

Now we have from definition of $L_n$ in equation \eqref{Ln} that
\[
\min_{1\le j\le n-m_n}\frac{S_{j}}{n+1-j}=\frac{\lambda_n}{n(1-\lambda_n)}L_n+\frac{k_n\lambda_n}{n},
\]
and thus obtain
\begin{eqnarray*}
\Delta_n&\le&\frac{\varepsilon_n}{1-\varepsilon_n}(\frac{\lambda_n}{n(1-\lambda_n)})^{-1}(\frac{\lambda_n}
{n(1-\lambda_n)}L_n+\frac{k_n\lambda_n}{n})\\
&=&\frac{\varepsilon_nL_n}{1-\varepsilon_n}+\frac{\varepsilon_n}{1-\varepsilon_n}k_n(1-\lambda_n)\\
&=&o_p(1)+O_p(\frac{\sqrt{k_n}(1-\lambda_n)}{\log n})\\
&=&o_p(1)
\end{eqnarray*}
from \eqref{1-lambda-upper}. Then it follows from \eqref{Deltan} that
$(\frac{\lambda_n}{n(1-\lambda_n)})^{-1}\Big(\max_{1\le j\le n-m_n}Y_{nj}^2-(1-\frac{k_n\lambda_n}{n})\Big)$ and $-L_n$ have the same asymptotic distribution. The lemma follows since $-L_n\xrightarrow{d}\Lambda$ from Lemma~\ref{maxSj}. \hfill $\blacksquare$

\vspace{10pt}

\noindent{\it Proof of Theorem~\ref{thm-main}.}
Set $\beta_n=1-\frac{k_n\lambda_n}{n}$ and $\alpha_n=\frac{\lambda_n}{n(1-\lambda_n)}$.  We apply Lemma~\ref{less} with $r_n=n-m_n$
and $p_n=n-k_n$ under condition~\eqref{cc}.  Note that $\lim_{n\to\infty}P(Y_{n1}^2>\beta_n+\alpha_n x)=0$ from Lemma~\ref{expansion1}.
 Then from Lemmas~\ref{expansion1} and \ref{maxY2} we have
\[
\frac{\max_{1\le j\le p_n}Y_{nj}^2-\beta_n}{\alpha_n}\xrightarrow{d} \Lambda.
\]
Then from Lemma~\ref{bird} we get
\[
\frac{\max_{1\le j\le n-m_n}Y_{nj}-\beta_n^{1/2}}{\alpha_n/(2\beta_n^{1/2})}\xrightarrow{d}\Lambda,
\]
i.e. Theorem ~\ref{thm-main} holds in view of \eqref{identical}.  This completes the proof.   \hfill $\blacksquare$

\vspace{20pt}

\noindent\textbf{Acknowledgements}. The authors would like to thank
the referee for his/her constructive suggestions for revision. The
research of Yu Miao was supported in part by NSFC (11971154). The
research of Yongcheng Qi was supported in part by NSF Grant
DMS-1916014.

\vspace{20pt}

%\noindent\textbf{Compliance with ethical standards}

%\noindent{\it Conflict of interest}  On behalf of all authors, the corresponding authors states that there is no conflict of interest.

\baselineskip 12pt
\def\ref{\par\noindent\hangindent 25pt}

\end{document}